\documentclass[12pt]{amsart}
\usepackage{amsmath,amsfonts,amstext,amsthm,amssymb,pxfonts}
\usepackage[colorlinks,citecolor=red,hypertexnames=false]{hyperref}
\usepackage{amsrefs}
\usepackage[pdftex]{graphicx}

\theoremstyle{plain} 
\newtheorem{lemma}[equation]{Lemma} 
 
\newtheorem{theorem}[equation]{Theorem} 
\newtheorem{corollary}[equation]{Corollary}

\theoremstyle{definition}

\theoremstyle{remark}
\newtheorem{remark}[equation]{Remark}

\allowdisplaybreaks

\numberwithin{equation}{section}

\title{Recovering Singular Integral Kernels from  Haar Shifts}

\author[ A.~Vagharshakyan ] {Armen Vagharshakyan }
\thanks{Research supported in part by NSF grant 0456611}
\address{ School of Mathematics, Georgia Institute of Technology, Atlanta GA 30332, USA}
\email{armenv@math.gatech.edu}
\begin{document}
\maketitle

\begin{center}
\textbf{October 22,2009}
\end{center}

\begin{abstract}
Any  sufficiently smooth one-dimensional Calder\'on-Zygmund  convolution operator is the average of 
Haar shift operators.  The latter are dyadic operators which can be efficiently expressed in terms 
of the Haar basis.  
 This extends the result of S. Petermichl \cite{MR1756958} on restoring Hilbert transform via Haar shift operators, 
 a technique that has become fundamental to the analysis of these operators. 
\end{abstract}
 
\section{Introduction.}

We will represent one-dimensional Calder\'on-Zygmund  convolution operators with sufficiently smooth kernels, see
\eqref{e.kernel_decrease}, \eqref{e.kernel_rapid_decrease},  by means of a properly chosen averaging of certain Haar shift operators with
bounded coefficients.  By \emph{Haar shift operators,} we mean for example \eqref{e.example}: linear operators 
that can be expressed in an efficient manner with the Haar basis. 

The use of Haar shift operators to represent singular integral operators goes back to the work of T.Figiel \cite{MR1110189}. Later
S.~Petermichl derived a strikingly succinct representation of the Hilbert transform \cite{MR1756958}. Similar representations were derived for Beurling \cite{MR1992955}, Riesz \cite{MR1964822} transforms and the truncated Hilbert transform (S. Petermichl, oral communication). The reason why these succinct representations are useful is that one can
deduce deep facts about singular integral operators, based on the analysis of Haar shift operators. 
In Petermichl's original paper \cite{MR1756958}, a deep property of Hankel operators associated to 
matrix symbols was deduced.  The linear $ A_2$ bound for the Hilbert transform was deduced 
by Petermichl \cite{MR2367098}, as was the same result for the Riesz transforms \cite{MR2367098} and the truncated Hilbert transform (S. Petermichl, oral communication).
The study of the Haar shift operators is interesting in itself 
\cites{MR1685781,MR2407233}, and has become an important model of the singular integral operators, 
see for instance their use in \cites{arXiv:0808.0832,MR2530853}.

To illustrate this, as a corollary to our main result, Theorem~\eqref{t.main} below, and 
the main result of \cite{arxiv:0906.1941}, we see that we have a proved 
a sharp $A_2$ inequality for the Calderon-Zygmund operators, a question of current interest: 
\begin{corollary}
Let 
\begin{equation*}
 T(f)(x)=P.V.\;\int_{\mathbb R} K(x-t) f(t) dt 
\end{equation*}
be a one dimensional Calderon-Zygmund convolution operator whose kernel  $K$ is odd and satisfies \eqref{e.kernel_decrease} and \eqref{e.kernel_rapid_decrease}, then 
\begin{equation*}
 \lVert Tf\rVert_{L_2(\omega)}\lesssim \lVert \omega \rVert_{A_2} \lVert f\rVert_{L_2(\omega)}.
\end{equation*}
By $ \lVert \omega \rVert_{A_2} $ we mean the $ A_2$ constant of the weight $ \omega $. 
(See \cites{MR2354322,MR2480568} for a definition.) 
\end{corollary}

This generalizes the result of S. Petermichl, obtained for the Hilbert transform \cite{MR2354322}, and improves the estimates of 
A. Lerner, S. Ombrosi and C. Perez \cite{MR2480568}*{equation 1.9} for these particular type of Calderon-Zygmund operators.

\section{Formulation of the Result.}
In order to formulate the main theorem we introduce some notations.

For any $\beta=\{\beta_l\}\in \{0,1\}^{\mathbb{Z}}$ and for any $r \in [1,2)$
define the dyadic grid ${\mathbb D}_{r,\beta}$ to be the collection of intervals 
\begin{equation*}
{\mathbb D}_{r,\beta}=
\Biggl\{
r 2^n
\left(
[0;1)+k+
\sum_{i<n}2^{i-n}\beta_i
\right)
\Biggr\}
_{n \in {\mathbb Z},\thinspace k\in {\mathbb Z}}
\end{equation*}
This parametrization of dyadic grids appears  explicitly in \cite{MR2464252}, and 
implicitly   in \cite{MR1998349}*{section 9.1}. 
Note, that the dyadic grid we use is different from the one used in \cite{MR1756958}. 

Place the usual uniform probability measure $\mathbb P $ on the space $\{0,1\}^{\mathbb{Z}}$, explicitly
\begin{equation*}
 \mathbb P (\beta: \beta_l=0)= \mathbb P (\beta: \beta_l=1)=\frac{1}{2}, \qquad \text{for all } l\in \mathbb{Z}.
\end{equation*}

We define two functions. Take $h$ to be the function supported on $[0,1]$ defined by
\begin{equation}\label{e.h}
h (x) = 
\begin{cases}
7,  & 0< x<1/4,
\\
-1,  & 1/4\le x \le 1/2, 
\\
1, & 1/2\le x\le 3/4,
\\
-7, & 3/4\le x\le 1. 
\end{cases}
\end{equation}
and  $g$ to be the function supported on $[0,1]$ defined by
\begin{equation} \label{e.g}
g (x) = 
\begin{cases}
-1, & 0\le x\le 1/4,
\\
1, & 1/4\le x\le 1/2,
\\
 1,  & 1/2\le x<3/4,
\\
-1,  & 3/4\le x \le 1. 
\end{cases}
\end{equation} 
Note that the function $g$ appears in \cite{MR1756958} paired with the usual Haar function.  In contrast, our function 
 $h$, defined by \eqref{e.h}, differs a little bit from the Haar function. In some sense, our choice of the function $h$ makes the convolution $h\ast g$
 'less smooth', and this property will be crucial for the proof. We'll make this statement precise in \eqref{e.a_is_invertible}, which will permit us to invert a Fourier transform.

For any function $f$ and any interval $I=[a;a+l]$ we define the function $f_I$ to be the scaling of $f$ to $I$ which preserves the $L_2$-norm, namely
\begin{equation}\label{e.shift}
f_I(x)=f_{[a,a+l]}=\frac{1}{\sqrt l}f\left(\frac{x-a}{l}\right).
\end{equation}

Now we are ready to state our main theorem:
\begin{theorem}\label{t.main}
Let $K:(-\infty,0)\cup(0,\infty)\rightarrow \mathbb R $ be an odd, twice differentiable function (in the sense that $K^{\prime}$ is absolutely continuous)
which satisfies
\begin{equation}\label{e.kernel_decrease}
 \lim_{x\rightarrow\infty}K(x)= \lim_{x\rightarrow\infty}K^{\prime}(x)=0 
\end{equation}
and
\begin{equation}\label{e.kernel_rapid_decrease}
x^3 K^{\prime\prime}(x)\in L_{\infty}(\mathbb R ).
\end{equation}
Then there exists a coefficient-function $\gamma \;:\; (0,\infty)\rightarrow  \mathbb R $, satisfying
\begin{equation*}
\lVert \gamma \rVert_{\infty} \leq C \lVert x ^3 K'''(x)\rVert_{\infty}
\end{equation*}
so that
\begin{equation}\label{e.kernel_restored}
K(x-y)=\int_{\{0,1\}^{\mathbb{Z}}}\int_{1}^2 \sum_{I\in {\mathbb D}_{r,\beta}}\gamma(|I|) \thinspace h_I (x) g_I (y)
\;  \frac{dr}{r} \thinspace d \mathbb P (\beta)
\end{equation}
for all $x\neq y$.  
Here $C$ is some absolute constant
and the series on the right of \eqref{e.kernel_restored}
is a.e. absolutely convergent.
\end{theorem}

\begin{remark}
Note, that for for fixed $ r, \beta $, and a function $ \gamma\in L_{\infty}(R) $, the linear operator 
\begin{equation} \label{e.example}
f \mapsto 
\sum_{I\in {\mathbb D}_{r,\beta}}\gamma(|I|) \langle g_I  , f  \rangle h_I (x)  
\end{equation}
is an example of a Haar shift operator as defined in \cite{arxiv:0906.1941}.  Note that this operator, expressed 
as a matrix in the Haar basis, has a bounded diagonals, but is even better than that: one only needs to 
use Haar coefficients associated with dyadic intervals that intersect and have lengths that differ 
by at most a factor of $ 2$. 
\end{remark}

\section{Proof of the Theorem} 

\subsection*{Step 1: Derivation of an Integral Equation.}
The following lemma derives a concise formula for properly averaged Haar shift operators:
\begin{lemma}\label{l.convolution}
Suppose the functions $h$ and $g$ are defined by \eqref{e.h},\eqref{e.g} and suppose $\gamma\in L_{\infty}(R_{+})$. Then for any $x\neq y$, we have
\begin{equation}\label{e.convolution}
\int_{\{0,1\}^{\mathbb{Z}}}\int_{1}^{2}   \sum_{I\in {\mathbb D}_{r,\beta}}\gamma(|I|) \thinspace h_I (x) g_I (y) \thinspace 
\frac{dr}{r} \thinspace d \mathbb P(\beta)=\int_{0}^{\infty}\frac{\gamma(r)}{r^2}\left(h\ast g_1\right)\left(\frac{x-y}{r}\right)dr
\end{equation}
where $g_1(x)\equiv g(-x)$, and the functions $h_I,g_I$ are defined by \eqref{e.shift}. 
Here the series on the left of \eqref{e.convolution} is a.e. absolutely convergent.
\end{lemma}
\begin{remark}
This lemma appears in \cite{MR2464252} for the case $\gamma\equiv 1$. 
\end{remark}
\begin{remark}
The notation $g_1$ is introduced in the right hand side of \eqref{e.convolution} to emphasize the role of convolution. 
\end{remark}

\begin{proof}
The following calculation justifies the a.e. convergence of series,
\begin{align*}
\sum_{I\in {\mathbb D}_{r,\beta}}\left|\gamma(|I|) \thinspace h_I (x) g_I (y) \right|
\leq 
\lVert \gamma \rVert_{\infty} \sum_{ \overset{I\in {\mathbb D}_{r,\beta}}{x\in I,\; y\in I}} \thinspace |h_I (x) g_I (y)|
\leq
\\
\leq 
\lVert \gamma \rVert_{\infty} \lVert h\rVert_{\infty}\lVert g\rVert_{\infty}\sum_{\overset{I\in {\mathbb D}_{r,\beta}}{ x\in I,\; |I|\geq |x-y|}} \frac{1}{|I|}
\leq 
\frac{2\lVert \gamma\rVert_{\infty}\lVert h\rVert_{\infty}\lVert g\rVert_{\infty}}{|x-y|}.
\end{align*}
Now recalling the definition of the dyadic grid $\mathbb D_{r,\beta}$ we get
\begin{align*}
\int_{\{0,1\}^{\mathbb{Z}}}\int_{1}^2 &
\sum_{I\in {\mathbb D}_{r,\beta}} \gamma(|I|) h_I (x) g_I (y)   \frac{dr}{r} d \mathbb P(\beta)
\\
& = 
\int_{\{0,1\}^{\mathbb{Z}}}\int_{1}^2 
\sum_{k \in {\mathbb Z}}\sum_{n \in {\mathbb Z}}
\frac{\gamma(r 2^{n})}{r 2^n}
h\left(\frac{x}{r 2^n}-k-\sum_{i<n}2^{i-n} \beta_i\right)
g\left(\frac{y}{r 2^n}-k-\sum_{i<n}2^{i-n} \beta_i\right)
 \frac{dr}{r} d \mathbb P(\beta)
\\
&= 
\int_{1}^2 
\sum_{n \in {\mathbb Z}}
\frac{\gamma(r 2^{n})}{r 2^n}
\int_{\mathbb R}
h\left(\frac{x}{r 2^n}-s\right)
g\left(\frac{y}{r 2^n}-s\right)
ds \frac{dr}{r}
\\&
=
\sum_{n \in {\mathbb Z}}\int_{1}^{2}
\frac{\gamma(r 2^{n})}{r 2^n}
\left(h \ast g_1\right)\left(\frac{x-y}{r 2^n}\right)\frac{dr}{r}
\\
&=
\sum_{n \in {\mathbb Z}}
\int_{2^n}^{2^{n+1}}
\frac{\gamma(r)}{r^2}
\left(h \ast g_1\right)\left(\frac{x-y}{r}\right)dr
\\&=
\int_{0}^{\infty}
\frac{\gamma(r)}{r^2}\left(h\ast g_1\right)\left(\frac{x-y}{r}\right)dr.
\end{align*}
\end{proof}

Having this lemma at hand, the claim of theorem \ref{t.main} is equivalent to the following:
For $h$ and $g$ defined by \eqref{e.h} and \eqref{e.g}, find a function $\gamma\in L_{\infty}(\mathbb R)$
which would satisfy the following integral equation
\begin{equation}\label{e.integralequation}
K(x)=
\int_0^{\infty} 
\frac{\gamma(r)}{r^2}
(h\ast g_1)\left(\frac{x}{r}\right)
dr,
\end{equation}
for all $x>0$. The case $x<0$ would be satisfied automatically as both $K$ and $h\ast g_1$ are odd.

\subsection*{     Step 2: Derivation of Recursive Equation}
In this step we'll use the functional equation \eqref{e.integralequation} to get a recursive equation \eqref{e.recursiveequation} for the coefficient function $\gamma$.

Differentiating \eqref{e.integralequation} twice, we get
\begin{equation*}
K^{\prime\prime}(x)=
\int_0^{\infty} 
\frac{\gamma(r)}{{r}^4}
\thinspace
(h\ast g_1)^{\prime\prime}\left(\frac{x}{r}\right)
dr, \qquad  x>0.
\end{equation*}
or equivalently
\begin{equation}\label{e.integralequationderived}
 x^3 K^{\prime\prime}(x)=\int_{0}^{\infty}t^2\gamma\left(\frac{x}{t}\right)(h\ast g_1)^{\prime\prime}(t)dt, \qquad x>0.
\end{equation}
Using the definitions \eqref{e.h},\eqref{e.g} for $h$ and $g$, we see that $h * g_1$ is 
a continuous, piecewise linear, odd function. The graph of $ h * g_1$  on the positive axis is illustrated in Figure 1.
\begin{figure}[ht]
\centering 
\includegraphics{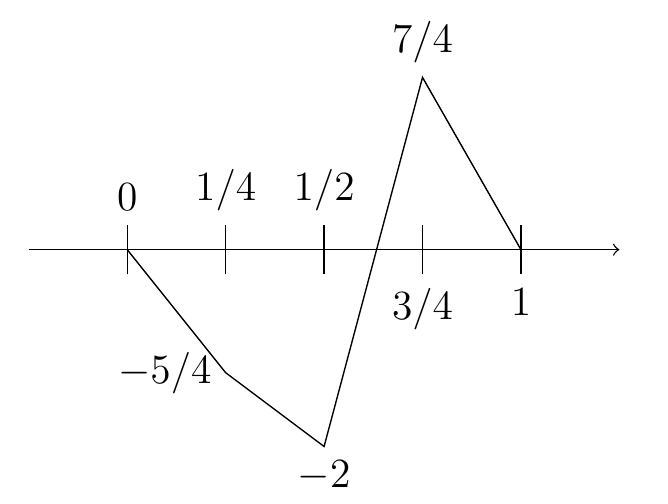}
\caption{$h\ast g_1$ on the positive axis.}
\end{figure}

  Thus, the  function $(h\ast g_1)^{\prime\prime}$ is a linear combination of Dirac measures, which one can calculate from the graph above, in particular
\begin{equation*}
 (h\ast g_1)^{\prime\prime}(x)=2\delta(x-1/4)+18\delta(x-1/2)-22\delta(x-3/4)+7\delta(x-1), \qquad x>0,
\end{equation*}
where $\delta$ is the usual Dirac delta function centered at the point $0$.

With this,  \eqref{e.integralequationderived} becomes
\begin{equation}\label{e.recursiveequation}
x^3 K^{\prime\prime}(x)=
2\left(\frac{1}{4}\right)^2\gamma(4x)
+18\left(\frac{1}{2}\right)^2\gamma(2x)
-22\left(\frac{3}{4}\right)^2\gamma\left(\frac{4}{3}x\right)
+7\gamma(x), \qquad  x>0.
\end{equation}
Let's modify \eqref{e.recursiveequation} to a form, which will be more convenient to us.
Denote
\begin{equation*}
m(x)=e^{3x}K^{\prime\prime}(e^x),
\end{equation*}
and
\begin{equation}\label{e.c}
c(x)=\gamma(e^x). 
\end{equation}
In terms of these new notations the equation \eqref{e.recursiveequation} becomes
\begin{equation}\label{e.c_recursively}
m(x)=
\tfrac{1}{8}c(x+\ln 4)
+\tfrac{9}{2} c(x+\ln 2)
-22\left(\tfrac{3}{4}\right)^2 c\left(x+\ln  \left(\tfrac43\right)\right)
+7 c(x), \qquad\ -\infty<x<\infty.
\end{equation}
The condition \eqref{e.kernel_rapid_decrease} of theorem \ref{t.main} provides that $m\in L_{\infty}(\mathbb R)$. We want to find $c\in L_{\infty}(\mathbb R)$ which would solve \eqref{e.c_recursively}.
\begin{remark}
In the case of Hilbert transform we have $m(x)\equiv 2$, thus a constant function  $c(x)\equiv C$ for a proper constant $C$ would solve  \eqref{e.c_recursively}.
\end{remark}

\subsection*{     Step 3: Fourier Transform}
We'll use Fourier transform in order to solve the recursive functional equation \eqref{e.c_recursively}. 
(Here we'll deal with Fourier transform of $L_{\infty}$ functions, which is understood in a distributional sense.)

Apply Fourier transform to both sides of \eqref{e.c_recursively} to get
\begin{equation}\label{e.fourierequation}
m^{\ast}(\omega)=a(\omega){c}^{\ast}(\omega),
\end{equation}
where 
\begin{equation*}
a(\omega)=
\tfrac{1}{8} e^{i\omega \ln 4}
+\tfrac{9}{2} e^{i\omega \ln  2}
-22\left(\frac{3}{4}\right)^2 e^{i\omega \ln  \left(\tfrac43\right)}
+7 .
\end{equation*}
Now the function $a$ is a Fourier transform of a finite Borel measure on $\mathbb R$.
Also note the following important property of $a$: our choice of functions $h$ and $g$ provided 
that one of the terms of $a$ dominates the rest
\begin{equation}\label{e.a_is_invertible}
12\tfrac38=22\left(\frac{3}{4}\right)^2 
>
\frac{1}{8}
+\frac{9}{2} 
+7=11\tfrac58.
\end{equation}
In particular, we have $ \lvert  a (\omega )\rvert\ge \tfrac 34 $ for all $ \omega $.

Recall that the space of Fourier transforms of finite Borel measures on $\mathbb R $, equipped with the 
$L_{\infty}$ norms of these Fourier transforms, is a Banach algebra under pointwise multiplication. 
Therefore $a$ is invertible, too. (The inverse of $a$ can be written in terms of a Neumann series of exponents.)
But this means that ${a} ^{-1} $ is a multiplier of the space $L_{\infty}(\mathbb R)$.
Hence, there exists a function $c\in L_{\infty}(\mathbb R)$, which solves the equation \eqref{e.fourierequation} and 
$\lVert c\rVert_{\infty}<C\lVert m\rVert_{\infty}$
(for some absolute constant $C$). 
Using \eqref{e.c} we can further restore the coefficient-function $\gamma$. It would solve the integral equation \eqref{e.integralequationderived} and would satisfy the same bound as the function $c$, i.e.
\begin{equation*}
\lVert\gamma\rVert_{\infty}\leq C\lVert m\rVert_{\infty}.
\end{equation*}
This fact, along with the conditions \eqref{e.kernel_decrease} on kernel $K$ is sufficient to make the integral in equation \eqref{e.integralequation} convergent, and to justify the passage from \eqref{e.integralequationderived} back to \eqref{e.integralequation}.

\begin{remark}\label{r.necessary}
The conditions \eqref{e.kernel_decrease} and \eqref{e.kernel_rapid_decrease} are somewhat necessary. Indeed, if some functions $h,g:[0,1]\rightarrow \mathbb R$ are constant on all dyadic intervals with sufficiently small length and if the coefficient function $\gamma$ is in $L_{\infty}(\mathbb R)$ then the lemma \eqref{l.convolution} still holds. Thus, whatever kernel $K$ is restored by the averaging of corresponding Haar shift operator, it must satisfy \eqref{e.integralequation} and \eqref{e.integralequationderived}.
If additionally $h$ is odd and $g$ is even with respect to the point $1/2$, then $h\ast g$ would be a piecewise linear function with bounded support, vanishing at $0$. So, \eqref{e.integralequation} and \eqref{e.integralequationderived} would imply that $K$ has to satisfy \eqref{e.kernel_decrease} and \eqref{e.kernel_rapid_decrease}.
\end{remark}


\begin{bibsection}
\begin{biblist}

\bib{MR1110189}{article}{
    AUTHOR = {Figiel, Tadeusz},
     TITLE = {Singular integral operators: a martingale approach},
 BOOKTITLE = {Geometry of {B}anach spaces ({S}trobl, 1989)},
    SERIES = {London Math. Soc. Lecture Note Ser.},
    VOLUME = {158},
     PAGES = {95--110},
 PUBLISHER = {Cambridge Univ. Press},
   ADDRESS = {Cambridge},
      YEAR = {1990},
}




\bib{MR2464252}{article}{
   author={Hyt{\"o}nen, Tuomas},
   title={On Petermichl's dyadic shift and the Hilbert transform},
   language={English, with English and French summaries},
   journal={C. R. Math. Acad. Sci. Paris},
   volume={346},
   date={2008},
   number={21-22},
   pages={1133--1136},
   issn={1631-073X},
   review={\MR{2464252}},
}

\bib{MR1992955}{article}{
    AUTHOR = {Dragi{\v{c}}evi{\'c}, Oliver}, 
    Author={Volberg, Alexander},
     TITLE = {Sharp estimate of the {A}hlfors-{B}eurling operator via
              averaging martingale transforms},
   JOURNAL = {Michigan Math. J.},
    VOLUME = {51},
      YEAR = {2003},
    NUMBER = {2},
     PAGES = {415--435},
      ISSN = {0026-2285},
}

\bib{arxiv:0906.1941}{article}{
author={Lacey, Michael},
author={Petermichl, Stephanie},
author={Reguera, Maria Carmen},
title={Sharp $A_2$ Inequality for Haar Shift Operators},
date={2009}, 
journal={Math. Ann., to appear}, 
eprint={http://arxiv.org/abs/0906.1941}
}

\bib{MR2530853}{article}{
   author={Lacey, Michael T.},
   author={Petermichl, Stefanie},
   author={Pipher, Jill C.},
   author={Wick, Brett D.},
   title={Multiparameter Riesz commutators},
   journal={Amer. J. Math.},
   volume={131},
   date={2009},
   number={3},
   pages={731--769},
   issn={0002-9327},
   review={\MR{2530853}},
}


     \bib{arXiv:0808.0832}{article}{
   author={Lacey, Michael T.},
   author={Pipher, Jill C.},
   author={Petermichl, Stefanie},
   author={Wick, Brett D.},
   title={Iterated Riesz Commutators: A Simple Proof of Boundedness},
   journal={to appear in Proceeding of El Escorial, 2008},
   eprint={http://www.arxiv.org/abs/0808.0832},
}   


\bib{MR2480568}{article}{
    AUTHOR = {Lerner, Andrei K.},
    AUTHOR={Ombrosi, Sheldy},
    AUTHOR={P{\'e}rez, Carlos},
     TITLE = {{$A\sb 1$} bounds for {C}alder\'on-{Z}ygmund operators related
              to a problem of {M}uckenhoupt and {W}heeden},
   JOURNAL = {Math. Res. Lett.},
    VOLUME = {16},
      YEAR = {2009},
    NUMBER = {1},
     PAGES = {149--156},
      ISSN = {1073-2780},
}


 \bib{MR1685781}{article}{
    author={Nazarov, F.},
    author={Treil, S.},
    author={Volberg, A.},
    title={The Bellman functions and two-weight inequalities for Haar
    multipliers},
    journal={J. Amer. Math. Soc.},
    volume={12},
    date={1999},
    number={4},
    pages={909--928},
    issn={0894-0347},
    review={\MR{1685781 (2000k:42009)}},
 }


\bib{MR2407233}{article}{
   author={Nazarov, F.},
   author={Treil, S.},
   author={Volberg, A.},
   title={Two weight inequalities for individual Haar multipliers and other
   well localized operators},
   journal={Math. Res. Lett.},
   volume={15},
   date={2008},
   number={3},
   pages={583--597},
   issn={1073-2780},
   review={\MR{2407233 (2009e:42031)}},
}

\bib{MR1998349}{article}{
   author={Nazarov, F.},
   author={Treil, S.},
   author={Volberg, A.},
   title={The $Tb$-theorem on non-homogeneous spaces},
   journal={Acta Math.},
   volume={190},
   date={2003},
   number={2},
   pages={151--239},
   issn={0001-5962},
   review={\MR{1998349 (2005d:30053)}},
}







\bib{MR1756958}{article}{
   author={Petermichl, Stefanie},
   title={Dyadic shifts and a logarithmic estimate for Hankel operators with
   matrix symbol},
   language={English, with English and French summaries},
   journal={C. R. Acad. Sci. Paris S\'er. I Math.},
   volume={330},
   date={2000},
   number={6},
   pages={455--460},
   issn={0764-4442},
   review={\MR{1756958 (2000m:42016)}},
}

\bib{MR2367098}{article}{
   author={Petermichl, Stefanie},
   title={The sharp weighted bound for the Riesz transforms},
   journal={Proc. Amer. Math. Soc.},
   volume={136},
   date={2008},
   number={4},
   pages={1237--1249},
   issn={0002-9939},
   review={\MR{2367098 (2009c:42034)}},
}

\bib{MR1964822}{article}{
   author={Petermichl, S.},
   author={Treil, S.},
   author={Volberg, A.},
   title={Why the Riesz transforms are averages of the dyadic shifts?},
   booktitle={Proceedings of the 6th International Conference on Harmonic
   Analysis and Partial Differential Equations (El Escorial, 2000)},
   journal={Publ. Mat.},
   date={2002},
   number={Vol. Extra},
   pages={209--228},
   issn={0214-1493},
   review={\MR{1964822 (2003m:42028)}},
}

\bib{MR2354322}{article}{
   author={Petermichl, S.},
   title={The sharp bound for the Hilbert transform on weighted Lebesgue
   spaces in terms of the classical $A\sb p$ characteristic},
   journal={Amer. J. Math.},
   volume={129},
   date={2007},
   number={5},
   pages={1355--1375},
   issn={0002-9327},
   review={\MR{2354322 (2008k:42066)}},
}






\end{biblist}
\end{bibsection}
\end{document}